\documentclass[12pt]{article}

\usepackage{fancyhdr}
\usepackage[lmargin=3cm,rmargin=3cm]{geometry}

\usepackage{graphicx}
\usepackage[brazil]{babel}
\usepackage[utf8]{inputenc}
\usepackage{amsmath,amssymb} 
\usepackage{amsthm}
\usepackage{amsfonts,amsmath,amscd}
\usepackage{xcolor}
\usepackage[colorlinks=true, bookmarksnumbered=true, bookmarksopen=true, bookmarksopenlevel=3, pdfstartview=FitH, linkcolor=blue, pdfmenubar=true, pdftoolbar=true, bookmarks=true,citecolor=red, urlcolor=blue, filecolor=magenta,plainpages=false,pdfpagelabels,breaklinks]{hyperref}

\newtheorem{theorem}{Theorem}[section]
\newtheorem{remark}[theorem]{Remark}

\newtheorem{corollary}[theorem]{Corollary}

\makeatletter
\renewcommand{\maketitle}{\vspace*{0pt}
	\begin{center}
		 \textbf{\@title}
	\end{center}
	\bgroup\setlength{\parindent}{0pt}
	\begin{flushright}
		\@author
	\end{flushright}\egroup
}

\newenvironment{abstracteng}{
	\begin{flushleft}%
		\bfseries{Abstract}
	\end{flushleft}}%

\numberwithin{equation}{section}

\begin{document}

\setcounter{page}{0}

\title{ {\bf \uppercase{Einstein Warped Products with Einstein Base and Fiber }}}
	\vspace{12pt}

\author{Márcio Lemes de Sousa \\{\small Universidade Federal de Mato Grosso, ICET-CUA, Pontal do Araguaia, Brasil,}\\{\small \underline{{\color{blue}marcio.sousa@ufmt.br} }}\\Tibério Bittencourt de Oliveira Martins \\{\small Universidade Federal de Mato Grosso, ICET-CUA, Pontal do Araguaia, Brasil,}\\{\small \underline{{\color{blue}tiberio.martins@ufmt.br} }}\\ Carlos Rodrigues da Silva \\{\small Universidade Federal de Mato Grosso, ICET-CUA, Pontal do Araguaia, Brasil,}\\{\small \underline{{\color{blue}carlos.silva@ufmt.br} }}}

\date{}

\maketitle
\noindent

	


\vspace{12pt}





	

\begin{abstracteng}
\noindent We study Einstein riemannian manifolds endowed with a warped product structure. We focus on the case in which both the base and the fiber are Einstein manifolds and establish necessary and sufficient conditions for the warped product itself to be an Einstein manifold. Moreover, we explicitly determine the warping function when the base is the hyperbolic space.
	
	\vspace{0.2cm}
	\end{abstracteng}

\section{Introduction}

A semi-riemannian manifold $(M^{n},g)$ is said to be Einstein if there exists a real constant $\lambda$ such that
$$Ric_{g}(X, Y) = \lambda g(X, Y)$$
for all $X, Y\in T_{p}M$ and all $p\in M$. This notion is relevant only for $n\geq 4$. In fact, if $n=1$, then $Ric_{g} = 0$. If $n=2$, then for each $p\in M$,
 $$Ric_{g}(X, Y) =
\frac{1}{2}Kg(X, Y),$$
so a 2-dimensional semi-riemannian manifold is Einstein if and only if it has constant sectional (or scalar) curvature.
According to do Carmo \cite[p. 83]{CAR}, for $n=3$, $(M^{3},g)$  is Einstein if and only if it has constant sectional curvature.

In recent years, several authors have considered the following problem:

Let $(M^{n}, g)$ be a semi-riemannian manifold of dimension $n>3$. Does there exist a metric $g'$ on $M$ such that $(M,g')$ is an Einstein manifold?

According to T. Aubin \cite{Au}, deciding whether a riemannian manifold admits an Einstein metric will be one of the important questions in riemannian geometry in the coming decades. Finding solutions to this problem is equivalent to solving a nonlinear system of second-order partial differential equations. In particular, semi-riemannian Einstein manifolds with vanishing Ricci curvature are special solutions, in the vacuum case ($T=0$), of the equation
\begin{equation}\label{eqeinst}
Ric_{g} - \frac{1}{2}K g = T
\end{equation}
where $K$ is the scalar curvature with respect to the metric $g$ and $T$ is a symmetric tensor of order $2$. If $g$ is a Lorentz metric on a 4-dimensional manifold, this is simply the Einstein field equation. When the tensor $T$ represents a physical field, such as an electromagnetic field, a perfect fluid, a pure radiation field, or the vacuum ($T=0$), the above equation has been studied in several works, most of which deal with solutions invariant under some symmetry group of the equation (see \cite{KRM} for details). When the metric $g$ is conformal to Minkowski spacetime, the solutions in the vacuum case are necessarily flat (see \cite{KRM}).

  Several authors have constructed new examples of Einstein manifolds. In \cite{ZIL}, Ziller constructed examples of compact manifolds with constant Ricci curvature. Chen, in \cite{CH}, constructed new examples of odd-dimensional Einstein manifolds, and in \cite{YAU}, Yau presented a survey on Ricci-flat manifolds. In \cite{Ku1}, Kühnel studied conformal transformations between Einstein spaces and, as a consequence of the results obtained, showed that there exists no riemannian Einstein manifold with non-constant sectional curvature that is locally conformally flat. This result was extended to semi-riemannian manifolds (see \cite{KR}). In this context, in order to construct examples of Einstein manifolds with non-constant sectional curvature, Sousa and Pina, in \cite{MR1}, worked with warped product manifolds that are not locally conformally flat.

Thus, let $(B,g_{B})$ and $(F,g_{F})$ be semi-riemannian manifolds and let $f>0$ be a smooth function on $B$. The warped product $M = B\times_{f}F$ is the product manifold $B\times F$ endowed with the metric tensor
$$\tilde{g} = g_{B} + f^{2}g_{F},$$
where $B$ is called the base of $M = B\times_{f}F$, $F$ is the fiber, and $f$ is the warping function. For example, polar coordinates determine a warped product structure in the case of spaces of constant curvature; in this case,
$$M=\mathbb{R}^{+}\times_{r}S^{n - 1} .$$

There are several studies relating warped product manifolds and locally conformally flat manifolds (see \cite{BGV1}, \cite{BGV2}, \cite{BGV3} and their references).

In several articles, the authors studied Einstein warped product manifolds under various curvature and symmetry conditions (see \cite{QC}, \cite{CPW} and \cite{CPW1}). In particular, He–Petersen–Wylie \cite{CPW} characterized riemannian Einstein warped product metrics when the base is locally conformally flat. 

In \cite{MR1}, the authors studied semi-riemannian Einstein manifolds with a warped product structure. They considered the case in which the base is conformal to a pseudo-euclidean space and invariant under the action of the $(n-1)$-dimensional translation group, while the fiber is Ricci-flat, and they classified all solutions in the Ricci-flat case. Moreover, they obtained explicit solutions in the vacuum case for the Einstein field equations.

In 2018, in \cite{LMR}, the authors analyzed semi-riemannian Einstein manifolds with a warped product structure $M = B^{n}\times _{f}F^{d}$ with $n\geq 3$ and $d\geq2$. When the base is compact and the fiber is Ricci-flat, they proved that $M$ reduces to a riemannian product. Furthermore, when the base is conformal to pseudo-euclidean space, has constant scalar curvature, and the fiber is a Ricci-flat semi-riemannian manifold, the authors determined all solutions in the case where the conformal factor is translation invariant.
It is well known that the Einstein condition in warped product geometries implies that the fiber must necessarily be an Einstein manifold (see \cite{Be}). 
In this work, we first classify Einstein warped product manifolds in the case where both the base and the fiber are Einstein manifolds (Theorem \ref{teor1}). We then present several consequences of this result, always under the assumption that both the base and the fiber are Einstein manifolds.
 
 First, we classify Ricci-flat manifolds when both the base and the fiber are Einstein manifolds (Corollary \ref{co1}). Under the same hypotheses, we classify warped products in the case where the base is compact (Corollary \ref{co2}). Finally, we classify Einstein warped products when the base is a complete Einstein manifold, assuming certain conditions on the difference between the Ricci curvatures of the base and of the warped product (Corollary \ref{co3}).

Considering $M=\mathbb{H}^{n}\times F^{d}$, where $\mathbb{H}^{n}$ denotes the $n$-dimensional hyperbolic space, we explicitly determine the warping function $f$ such that the warped product $M$ is an Einstein manifold (Theorem \ref{teor2}).

Finally, under the hypotheses of Theorem \ref{teor2}, we establish the necessary conditions for the warping function $f$ to be globally defined and conclude that the Ricci curvature of the fiber must satisfy $\lambda_{F}\leq 0$ (Corollaries \ref{co4}  and \ref{co5}).

\begin{remark}
When the dimension of the fiber $F^{d}$ is $d = 1$, we consider $$M =
(\mathbb{R}^n, \overline{g})\times _{f}\mathbb{R},$$ and, in this case,
$\lambda_{F} = 0$.
\end{remark}

In what follows, we state our main results. We denote by
$\varphi_{ ,x_ix_j}$ and $f_{ ,x_ix_j}$ the second-order derivatives of $\varphi$ and $f$ with respect to $x_ix_j$.

\section{Main Results}
\begin{theorem}\label{teor1}
Let $(B^{n}, g_{B})$ and $(F^{d}, g_{F})$, with $n \geq 3$, be Einstein manifolds whose Ricci curvatures are constant and equal to $\lambda_{B}$ and $\lambda_{F}$, respectively. Then, the warped product manifold $M = B \times_{f} F$, with $f$ non-constant, is an Einstein manifold with constant Ricci curvature $\lambda$ if and only if the following conditions $(i)$, $(ii)$, $(iii)$ and $(iv)$ are satisfied:
\begin{itemize}
  \item[(i)] $Hess_{g_{B}}f(X,Y)= \displaystyle\frac{f}{d}(\lambda_{B} - \lambda)g_{B}(X, Y),\forall\ X, Y\in{\cal L}(B)$;

  \item[(ii)] $\lambda = (1 + \frac{d}{n-1})\lambda_{B}$;

  \item[(iii)]The function $f$ satisfies
   $$\|\nabla f\|^{2} = -\rho f^{2} + c, $$
  where $c$ is a constant and $\rho$ is the normalized scalar curvature of $B$.
  \item[(iv)] $\lambda_{F} = (d - 1)c.$

\end{itemize}
\end{theorem}

\begin{corollary}\label{co1}
Let $B$, $F$ and $f$ be as in the hypotheses of Theorem \ref{teor1}.
Then, the warped product manifold $M = B \times_{f} F$ is Ricci-flat if and only if the base manifold $B$ is also Ricci-flat.
\end{corollary}

\begin{corollary}\label{co2}
  Let $B$, $F$ and $f$ be as in the hypotheses of Theorem \ref{teor1}.
If $B$ is compact and the warped product $M = B \times_{f} F$ is an Einstein manifold, then $(B, g_{B})$ is isometric to a sphere of some radius.
\end{corollary}

\begin{corollary}\label{co3}
Let $B$, $F$ and $f$ be as in the hypotheses of Theorem \ref{teor1}, with $B$ complete.
If the warped product $M = B \times_{f} F$ is an Einstein manifold and
$$\frac{\lambda_{B} - \lambda}{d} = -r^{2},$$
then $(B, g_{B})$ is isometric to a sphere of radius $1/r$.
\end{corollary}

In the next theorem, we will use the hyperbolic space, defined by
$$\mathbb{H}^{n} = \{x = (x_{1}, x_{2}, \ldots, x_{n-1}, x_{n})\in\mathbb{R}^{n} ; x_{n}>0\}$$

\begin{theorem}\label{teor2}
  Let $(\mathbb{H}^{n}, g_{-1})$ be the $n$-dimensional hyperbolic space, with $n \geq 3$, whose metric is given by
$$(g_{-1})_{ij}=\frac{\delta{ij}}{x_{n}^{2}}.$$
Consider the warped product manifold
$$M = (\mathbb{H}^{n}, g_{-1})\times_{f}F^{d},$$
where $F^{d}$ is a riemannian Einstein manifold with constant Ricci curvature $\lambda_{F}$.
Then $M$ is an Einstein manifold with constant Ricci curvature $\lambda < 0$ if and only if
$$f = \frac{1}{x_{n}}\left(\sum_{j = 1}^{n -1}\frac{a}{2}x_{j}^{2} + b_{j}x_{j} + c_{j}\right) + \frac{a}{2}x_{n} + \frac{b}{x_{n}},$$
where $a$, $b$, $b_{j}$ and $c_{j}$, for $j = 1, \ldots, n - 1$, are constants.
\end{theorem}

\begin{corollary}\label{co4}
Under the hypotheses of Theorem \ref{teor2}, if $a > 0$, $b \geq 0$ and
$$b_{j}^{2} - 2ac_{j} \leq 0,\ \ \ \hbox{for all}\ j = 1, \ldots, n - 1,$$
then $f$ is globally defined and $\lambda_{F} \leq 0$.
\end{corollary}

\begin{corollary}\label{co5}
In the context of Theorem \ref{teor2}, if $a = 0$ and $b \geq 0$, with
$$\sum_{j = 1}^{n - 1} c_{j} + b > 0\ \ \hbox{e}\ \ b_{j} = 0,$$
then $f$ is globally defined and $F$ is a Ricci-flat manifold.
\end{corollary}

\section{Proofs of the Main Results}

\textbf{Proof of the Theorem \ref{teor1}:} It follows from \cite{O'neil} that, if $X, Y \in \mathcal{L}(B)$ and $Z, W \in \mathcal{L}(F)$, then
\begin{eqnarray}
\label{ric1}
  Ric_{g}(X,Y) &=& Ric_{g_{B}}(X,Y) - \frac{d}{f}Hess_{g_{B}}f(X, Y) \nonumber\\
    Ric_{g}(X, Z) &=& 0  \\
    Ric_{g}(Z, W) &=& Ric_{g_{F}}(Z, W) -
 g(Z, W)\left(\frac{\triangle_{g_{B}}f}{f} + (d-1)\frac{g_{B}(\nabla f, \nabla
  f)}{f^{2}}\right). \nonumber
\end{eqnarray}
Assume initially that $(M,g)$ is Einstein. In this case, the first equation in \eqref{ric1} yields 
 $$Ric_{g_{B}}(X,Y) - \frac{d}{f}Hess_{g_{B}}f(X, Y) = \lambda g(X, Y).$$
 Since, by hypothesis, $B$ is Einstein and $g(X, Y) = g_{B}(X, Y)$, we obtain
 \begin{equation}\label{hess}
  Hess_{g_{B}}f(X,Y) = \displaystyle\frac{f}{d}(\lambda_{B} - \lambda)g_{B}(X, Y).
\end{equation}
Thus, the Hessian of $f$ with respect to the metric $g_{B}$ is a multiple of the metric itself. Moreover, since $B$ is Einstein by hypothesis, it follows from \cite{Ku1} that
 \begin{equation}\label{lap}
  \frac{f}{d}(\lambda_{B} - \lambda) = \frac{\triangle_{g_{B}} f}{n} = -\rho f.
  \end{equation}
Consequently, $$\displaystyle \rho = \frac{(\lambda - \lambda_{B})}{d}.$$
On the other hand, since $B$ is Einstein, we also have
$$\displaystyle\rho = \frac{\lambda_{B}}{n - 1}.$$
Therefore,
\begin{equation}\label{cric}
  \lambda = (1 + \frac{d}{n - 1})\lambda_{B}.
\end{equation}
Again, by \cite{Ku1}, we have 
\begin{equation}\label{grad}
\|\nabla f\|^{2} = -\rho f^{2} + 2bf + c.
\end{equation}
where $b$ and $c$ are constants. We shall show that, in this case, $b=0$.

Since $(F, g_{F})$ is Einstein,  $g(Z, W) = f^{2}g_{F}(Z, W)$  and $M$ is also Einstein, the third equation in \eqref{ric1} gives $$\lambda_{F}g_{F}(Z, W) - \left( \frac{\triangle_{g_{B}}f}{f} +(d - 1)\frac{\|\nabla f\|^{2}}{f^{2}}\right)f^{2}g_{F}(Z, W) = \lambda f^{2}g_{F}(Z, W),$$
  that is,
\begin{equation}\label{fibra}
\lambda_{F} -  f\triangle_{g_{B}}f - (d - 1)\| \nabla f\|^{2} = \lambda f^{2}.
\end{equation}
Substituting \eqref{lap} and \eqref{grad} into \eqref{fibra}, we obtain
$$\lambda_{F} - f^{2}\frac{n}{d}(\lambda_{B} - \lambda) - (d - 1)(-\rho f^{2}+ 2bf + c) = \lambda f^{2}.$$
Thus, we arrive at 
\begin{equation}\label{fibra1}
\left[\lambda_{F} - (d- 1)c\right] - 2bf + \left[\frac{n}{d}(\lambda - \lambda_{B}) + \rho(d - 1)\right]f^{2} = \lambda f^{2}.
\end{equation}
Since $\rho=\dfrac{\lambda_{B}}{n-1}$, the coefficient of $f^{2}$ in \eqref{fibra1} can be rewritten as
$$\frac{n}{d}(\lambda - \lambda_{B}) + \rho(d - 1) = \frac{n}{d}(\lambda - \lambda_{B}) + (d - 1)\frac{\lambda -  \lambda_{B}}{d} = \frac{1}{d}(n + d -1)(\lambda - \lambda_{B}).$$
Solving \eqref{cric} for  $\lambda_{B}$ and substituting it into the right-hand side of the equality above, we obtain 
\begin{equation}\label{base}
\frac{n}{d}(\lambda - \lambda_{B}) + \rho(d - 1) = \frac{1}{d}(n + d -1)(\lambda - \frac{n-1}{n +d -1}\lambda) =\lambda.
\end{equation}
Thus, it follows from \eqref{base} that equation \eqref{fibra1} is satisfied if and only if $b = 0$ and
\begin{equation}\label{fric}
\lambda_{F} = (d- 1)c.
\end{equation}
For the converse, it suffices to recall that $B$ and $F$  are Einstein manifolds and to substitute expressions \eqref{hess}, \eqref{lap}, \eqref{grad}, and \eqref{fric} into the first and third equations of \eqref{ric1}. This leads directly to the desired result.

Therefore, the proof of Theorem \ref{teor1} is complete.

\hfill $\Box$

\textbf{Proof of the  Corollary \ref{co1}:} A straightforward consequence of condition $(ii)$ of Theorem \ref{teor1} is that $\lambda=0$ if and only if $\lambda_{B}=0$.

\hfill $\Box$

\textbf{Proof of the  Corollary \ref{co2}:}
By hypothesis, $B$ is a compact Einstein manifold; in particular, $B$ has constant scalar curvature. Then, by \cite{Ku1}, since the warped product $M= B\times_{f}F$ is Einstein, it follows that $B$ is a sphere of some radius.

\hfill $\Box$

\textbf{Proof of the  Corollary \ref{co3}:} By hypothesis, $B$ is a complete Einstein manifold. If the warped product $M= B\times_{f}F$ is Einstein with $$\displaystyle\frac{\lambda_{B} - \lambda}{d} = -r^{2},$$ then
$$Hess_{g_{B}}f(X,Y) = -r^{2}fg_{B}(X, Y)$$
for any tangent vectors $X$ and $Y$ on $B$. It follows from \cite{Ku1} that $B$ is a sphere of radius $1/r$.

\hfill $\Box$

\textbf{Proof of the  Theorem \ref{teor2}:}
Assume initially that $M$ is Einstein. Since $\mathbb{H}^{n}$ is Einstein with Ricci curvature $\lambda_{\mathbb{H}^{n}} = -(n - 1)$, by Theorem \ref{teor1} we have
$$\lambda = -(n + d -1) < 0$$
and
$$Hess_{g_{-1}}f(X, Y) = f\cdot g_{-1}(X, Y).$$
Let $\left\{X_{1}, X_{2} \ldots X_{n-1}, X_{n}\right\}$ be a basis of $\mathbb{H}^{n}$ such that $$g_{-1}(X_{i}, X_{j}) = \frac{\delta_{ij}}{x_{n}^{2}}.$$
Then,
\begin{equation}\label{hess2}
Hess_{g_{-1}}f(X_{i}, X_{j}) = f\frac{\delta_{ij}}{x_{n}^{2}}.
\end{equation}
On the other hand, according to \cite{O'neil}, we have 
\[
Hess_{g_{-1}}f(X_{i},X_{j})=f_{ ,x_ix_j}-\sum_k
\Gamma_{ij}^k f_{,x_k},
\]
where $\Gamma_{ij}^k$ are the Christoffel symbols of the metric $g_{-1}$. For $i,\ j,\ k$ pairwise distinct, since $(\mathbb{H}^{n}, g_{-1})$ is locally conformally flat, it follows from \cite{CAR} that
\[
\Gamma_{ij}^k= 0\ \ \ \ \ \ \ \ \ \Gamma_{ij}^i= -\frac{\varphi_{,x_{j}}}{\varphi}\ \ \ \ \ \ \ \ \ \ \Gamma_{ii}^k= \frac{\varphi_{,x_{k}}}{\varphi}\ \ \ \ \ \ \ \ \ \Gamma_{ii}^i= -\frac{\varphi_{,x_{i}}}{\varphi}. \] 
Therefore,
\begin{equation}\label{hes1}
\left\{
\begin{array}{ccl}
Hess_{g_{-1}}f(X_{i},X_{j}) &=& \displaystyle
f_{,x_ix_j}+\frac{\varphi_{,x_j}}{\varphi}f_{,x_i}
+\frac{\varphi_{,x_i}}{\varphi}f_{,x_j},\ \forall\ i\neq j = 1\ldots n\\
Hess_{g_{-1}}f(X_{i},X_{i}) &=& \displaystyle
f_{,x_ix_i}+2\frac{\varphi_{,x_i}}{\varphi}f_{,x_i} -
\sum_{k=1}^{n}\frac{\varphi_{,x_k}}{\varphi}f_{,x_k}.
\end{array}
\right.
\end{equation}
In the case $\varphi = x_{n}$, we obtain
\begin{equation}\label{hes2}
\left\{
\begin{array}{ccl}
Hess_{g_{-1}}f(X_{i},X_{j}) &=&
f_{,x_ix_j},\ \forall\ i\neq j \neq n \neq i\\
Hess_{g_{-1}}f(X_{i},X_{n}) &=& f_{,x_ix_n}+\frac{1}{x_{n}}f_{,x_i},\ \forall\ i\neq  n\\
Hess_{g_{-1}}(f)(X_{i},X_{i}) &=&
f_{,x_ix_i} -\frac{1}{x_{n}}f_{,x_n}\ \forall\ i\neq  n\\
Hess_{g_{-1}}(f)(X_{n},X_{n}) &=&
f_{,x_nx_n}+ \frac{1}{x_{n}}f_{,x_n}.
\end{array}
\right.
\end{equation}
Therefore, substituting \eqref{hes2} into \eqref{hess2}, we obtain 
\begin{equation}\label{hess3}
\left\{
\begin{array}{rcl}
f_{,x_ix_j} &=& 0,\ \forall\ i\neq j \neq n \neq i\\
x_{n}f_{,x_ix_n}+ f_{,x_i} &=& 0,\ \forall\ i\neq  n\\
x_{n}^{2}f_{,x_ix_i} - x_{n}f_{,x_n}&=& f,\forall\ i\neq  n\\
x_{n}^{2}f_{,x_nx_n} + x_{n}f_{,x_n} &=& f.
\end{array}
\right.
\end{equation}
Integrating the second equation in \eqref{hess3}, we obtain 
\begin{equation}\label{ed}
f_{,x_{i}} = \frac{l_{i}(\hat{x}_{n})}{x_{n}},\  \forall\ i\neq  n,
\end{equation}
where $l_{i}(\hat{x}_{n})$ is an arbitrary function independent of the variable $x_{n}$. Differentiating \eqref{ed} with respect to $x_{j}$, with $i,j, n$ pairwise distinct, and using the first equation in \eqref{hess3}, we conclude that $l_{i}(\hat{x}_{n}) = l_{i}(x_{i})$. Hence, integrating \eqref{ed} with respect to $x_{i}$, we obtain
\begin{equation}\label{ed2}
f = \frac{1}{x_{n}}\int l_{i}dx_{i} + m(\hat{x}_{i}).
\end{equation}
where $m(\hat{x}_{i})$ is an arbitrary function independent of the variable $x_{i}$. Fixing $i$ and differentiating \eqref{ed2} with respect to $x_{j}$, with $i, j, n$ pairwise distinct, and using \eqref{ed}, we have 
$$ m_{,x_{j}} = \frac{l_{j}(x_{j})}{x_{n}}$$
that is,
$$f = \frac{1}{x_{n}}\int l_{i}dx_{i} + \frac{1}{x_{n}}\int l_{j}dx_{j} +m_{ij}(\hat{x}_{i}, \hat{x}_{j})$$
where $m_{ij}(\hat{x}_{i}, \hat{x}_{j})$ is an arbitrary function independent of the variables $x_{i}$ e $x_{j}$. Proceeding inductively, we conclude that
\begin{equation}\label{torc}
f = \frac{1}{x_{n}}\sum_{i = 1}^{n -1}\int l_{i}dx_{i} + r({x_{n}}),
\end{equation}
where $r$ is an arbitrary function of $x_{n}$. 

Computing the partial derivative of $f_{,x_{i}}$ with respect to $x_{i}$ and the partial derivative of $f$ with respect to $x_{n}$, we obtain
\begin{equation}\label{part}
f_{,x_ix_i} = \frac{l_{i}'}{x_{n}},\ \forall\ \ i\neq n\ e\ \  f_{,x_{n}} = -\frac{1}{x_{n}}(f - r) + r'.
\end{equation}
Substituting \eqref{part} into the third equation of \eqref{hess3}, we obtain 
$$x_{n}^{2} \frac{l_{i}'}{x_{n}} - x_{n}(-\frac{1}{x_{n}}(f - r) + r')= f,$$
that is, 
$$l_{i}' = r' + \frac{1}{x_{n}}r. $$
Since $l_{i} = l_{i}(x_{i})$ and $r = r(x_{n})$, with $i\neq n$, there exists a constant $a$ such that
\begin{equation}\label{ed3}
l_{i}' = a\ e\ \  r' + \frac{1}{x_{n}}r = a.
\end{equation}

Integrating both equations in \eqref{ed3}, we obtain
\begin{equation}\label{ed4}
l_{i}(x_{i}) = ax_{i} + b_{i},\ \forall\ i\neq n, \ \hbox{e}\ \ r(x_{n}) = \frac{ax_{n}}{2} + \frac{b}{x_{n}},
\end{equation}
where $b_{i}$ and $b$ are constants.

Substituting the functions $l_{i}(x_{i})$ e $r(x_{n})$, obtained in \eqref{ed4}, into the expression for $f$, in \eqref{torc}, we arrive at 
\begin{equation}\label{tor}
f = \frac{1}{x_{n}}\left(\sum_{j = 1}^{n -1}\frac{a}{2}x_{j}^{2} + b_{j}x_{j} + c_{j}\right) + \frac{a}{2}x_{n} + \frac{b}{x_{n}},
\end{equation}
 where $a$, $b$, $b_{j}$, $c_{j}$, $ j = 1,\ldots, n - 1$, are constants. With this expression, the fourth equation in \eqref{hess3} is trivially satisfied.

Conversely, computing $Hess_{g_{-1}}f(X_{i},X_{j})$ from \eqref{tor} for $i,j=1,\ldots,n$,  one verifies that \eqref{hess2} is satisfied.

This concludes the proof of Theorem \ref{teor2}.

 \hfill $\Box$
 
 \textbf{Proof of the  Corollary \ref{co4}:} Since $x_{n} >0$, $a >0$ and $b\geq 0$, it follows that
\begin{equation}\label{f1}
 \frac{ax_{n}}{2} + \frac{b}{x_{n}} >0.
\end{equation}

On the other hand, since $a>0$ and $b_{j}^{2} - 2ac_{j} \leq 0 $, $\forall\ j = 1,\ldots, n-1 $, we have
\begin{equation}\label{f2}
\frac{a}{2}x_{j}^{2} + b_{j}x_{j} + c_{j}\geq 0.
\end{equation}
Combining \eqref{f1} and \eqref{f2}, we conclude that $f$ is positive on all of $\mathbb{H}^{n}$, hence globally defined.

Now, in order to compute $\lambda_{F}$, recall that
$$f_{,x_{i}} = \frac{ax_{i} + b_{i}}{x_{n}},\ \forall\ i \neq n,\ \ \hbox{and}\ \ f_{,x_{n}} = -\frac{1}{x_{n}}f + a.$$
 According to \cite{O'neil}, the gradient of the function $f$ is given by
$$\nabla f=\sum_{i,j=1}^{n}g^{ij}f_{,x_{i}} \frac{\partial }{\partial x_{j}},$$
where $(g^{ij})$ denotes the inverse metric and $\frac{\partial }{\partial x_{j}}$ are elements of a basis of the tangent space. Therefore,
\begin{equation}\label{grad2}
\| \nabla f\|^{2} = x_{n}^{2}\sum_{j = 1}^{n}f _{, x_{j}}^{2}.
\end{equation}
 Hence,
$$\|\nabla f\|^{2} =  x_{n}^{2}\sum_{j = 1}^{n - 1}\left(\frac{ax_{j} + b_{j}}{x_{n}}\right)^{2} + x_{n}^{2}\left(-\frac{1}{x_{n}}f + a\right)^{2},$$
that is,
$$ \|\nabla f\|^{2} = 2a\sum_{j = 1}^{n - 1}\left(\frac{a}{2}x_{j}^{2} + b_{j}\right) + \sum_{j = 1}^{n - 1}b_{j}^{2} + f^{2} -2ax_{n}f + a^{2}x_{n}^{2}.$$
Since $$\sum_{j = 1}^{n - 1}\left(\frac{a}{2}x_{j}^{2} + b_{j}\right) = x_{n}f - \frac{a}{2}x_{n}^{2}- b - \sum_{j=1}^{n-1}c_{j},$$ we obtain
$$\|\nabla f\|^{2} = 2a(x_{n}f - \frac{a}{2}x_{n}^{2}- b - \sum_{j=1}^{n-1}c_{j}) + \sum_{j = 1}^{n - 1}b_{j}^{2} + f^{2} -2ax_{n}f + a^{2}x_{n}^{2},$$
that is,
$$\|\nabla f\|^{2} = f^{2} + \left(\sum_{j = 1}^{n-1}(b_{j}^{2} -2ac_{j}) - b\right).$$
It follows from Theorem \ref{teor1} that $\lambda_{F} = (d -1)c$, where $$c =\displaystyle \sum_{j = 1}^{n-1}(b_{j}^{2} -2ac_{j}) - b.$$ Since $d - 1 \geq 0$, , we have $\lambda_{F} \leq 0$, because $(b_{j}^{2} -2ac_{j})\leq 0$ and $b\geq 0$.
This concludes the proof of Corollary \ref{co4}.

 \hfill $\Box$

\textbf{Proof of the  Corollary \ref{co5}} Under the hypotheses of this corollary, it is easy to see that $f$ is positive, and hence globally defined. On the other hand,
$$f_{,x_{i}} = 0\ \forall\ i \neq n\ \ \hbox{and}\ \ f_{,x_{n}} = -\frac{1}{x_{n}}f.$$
Then, by equation \eqref{grad2}, we have 
$$\|\nabla f\|^{2} = f^{2}.$$
Therefore, by the same argument as in Corollary \ref{co4}, it follows that $\lambda_{F} = 0$. This concludes the proof of Corollary \ref{co5}.

 \hfill $\Box$


\end{document}